\documentclass[11pt]{article}
\usepackage{amssymb}
\usepackage{amsmath}
\usepackage{graphicx}
\usepackage{amsfonts}
\usepackage{theorem}

\newtheorem{Lem}{Lemma}[section]

\newtheorem{The}[Lem]{Theorem}
\newtheorem{Prop}[Lem]{Proposition}
\newtheorem{Cor}[Lem]{Corollary}

\newtheorem{Rem}[Lem]{Remark}

\input{tcilatex}
\begin{document}

\title{Fej\'{e}r-type inequalities }
\author{Nicu\c{s}or Minculete$^{1}$\thanks{%
E-mail:minculeten@yahoo.com\newline
} \ and Flavia-Corina Mitroi$^{2}$\thanks{%
E-mail:fcmitroi@yahoo.com} \\
$^{1}${\small "Dimitrie Cantemir" University, }\\
{\small 107 Bisericii Rom\^{a}ne Street, Bra\c{s}ov, 500068, Romania}\\
$^{2}${\small University of Craiova, Department of Mathematics,}\\
{\small Street A. I. Cuza 13, Craiova, RO-200585, Romania}}
\date{}
\maketitle

\textbf{Abstract.} The aim of this paper is to present some new Fej\'{e}%
r-type results for convex functions. Improvements of Young's inequality (the
arithmetic-geometric mean inequality) and other applications to special
means are pointed as well.

\textbf{Keywords : } Fej\'{e}r's inequality, convex function,
Hermite-Hadamard inequality, Young's inequality.

\textbf{2010 Mathematics Subject Classification : } 26A51, 26D15

\section{Preliminaries}

{We start from an important result related to the convex functions due to
Ch. Hermite \cite{Her1883} and J. Hadamard \cite{Had1893} which asserts that
for every continuous convex func}tion{\ \textit{\ }}$f:[a,b]\rightarrow 
\mathbb{R}${\textit{\ }}the following inequalities hold:{\ 
\begin{equation}
f\left( \frac{a+b}{2}\right) \leq \frac{1}{b-a}\int_{a}^{b}f(x)\mathrm{d}%
x\leq \frac{f(a)+f(b)}{2}.
\end{equation}%
}

Fej\'{e}r \cite{Fej1906} established the following well-known weighted
generalization:

\begin{Prop}
\label{prop_fej}If {\textit{\ }}$f:[a,b]\rightarrow \mathbb{R}${\textit{\ \
is }continuous and convex and if }$g:[a,b]\rightarrow \mathbb{R}_{+}$\textit{%
\ }{\ is integrable and symmetric about }$\frac{a+b}{2}$ ( i.e. $g\left(
x\right) =g\left( a+b-x\right) $)$,$ then{\ }the following inequalities hold:%
{%
\begin{equation}
f\left( \frac{a+b}{2}\right) \int_{a}^{b}g(x)\mathrm{d}x\leq
\int_{a}^{b}f(x)g\left( x\right) \mathrm{d}x\leq \frac{f(a)+f(b)}{2}%
\int_{a}^{b}g(x)\mathrm{d}x.  \label{fej}
\end{equation}%
}
\end{Prop}

Before stating the results we recall some useful facts from the literature.
S. S. Dragomir, P. Cerone and A. Sofo present in \cite{Dra2000a,Dra2000b}
the following estimates of the precision in the Hermite-Hadamard inequality:

\begin{Prop}
\label{nic}Let $f:[a,b]\rightarrow \mathbb{R}$ be a twice differentiable
function such that there exists real constants $m$ and $M$ so that $m\leq
f^{\prime \prime }\leq M$. Then%
\begin{equation}
m\frac{(b-a)^{2}}{24}\leq \frac{1}{b-a}\int_{a}^{b}f(x)\mathrm{d}x-f\left( 
\frac{a+b}{2}\right) \leq M\frac{(b-a)^{2}}{24}  \label{cpn_1}
\end{equation}%
\textit{\ and}%
\begin{equation}
m\frac{(b-a)^{2}}{12}\leq \frac{f(a)+f(b)}{2}-\frac{1}{b-a}\int_{a}^{b}f(x)%
\mathrm{d}x\leq M\frac{(b-a)^{2}}{12}.  \label{cpn_2}
\end{equation}%
\textit{\ }
\end{Prop}

These inequalities follow from the Hermite-Hadamard inequality, for the
convex functions $f(x)-m\frac{x^{2}}{2}$ and $M\frac{x^{2}}{2}-f\left(
x\right) $.\textit{\ }

Motivated by the above results, the purpose of this paper is to discuss
further inequalities of Fej\'{e}r type.

\section{Fej\'{e}r type inequalities for convex functions}

\begin{The}
\label{th_2}Let $f:[a,b]\rightarrow \mathbb{R}$ be a twice differentiable
function such that there exist real constants $m$ and $M$ so that $m\leq
f^{\prime \prime }\leq M$. Then 
\begin{equation}
m\frac{\lambda (1-\lambda )}{2}(a-b)^{2}\leq \lambda f(a)+(1-\lambda
)f(b)-f(\lambda a+(1-\lambda )b)\leq M\frac{\lambda (1-\lambda )}{2}%
(a-b)^{2},  \label{ref_1}
\end{equation}%
for all $\lambda \in \lbrack 0,1].$
\end{The}

\noindent \textit{Proof:} We consider the function $g:[0,1]\rightarrow 
\mathbb{R}$, defined by 
\begin{equation*}
g(\lambda )=\lambda f(a)+(1-\lambda )f(b)-f(\lambda a+(1-\lambda )b)-m\frac{%
\lambda (1-\lambda )}{2}(a-b)^{2}.
\end{equation*}%
Since 
\begin{equation*}
g^{\prime \prime }(\lambda )=(a-b)^{2}[m-f^{\prime \prime }(\lambda
a+(1-\lambda )b)]\leq 0,
\end{equation*}%
the function $g$ is concave. But $g(0)=g(1)=0$, which implies that 
\begin{equation*}
0=(1-\lambda )g(0)+\lambda g(1)\leq g((1-\lambda )\cdot 0+\lambda \cdot
1)=g(\lambda ),
\end{equation*}%
for all $\lambda \in \lbrack 0,1]$. Therefore, we obtain the first part of
inequality (\ref{ref_1}).\newline
To see that the later inequality holds, our next step is to take the convex
function $h:[0,1]\rightarrow \mathbb{R}$, defined by 
\begin{equation*}
h(\lambda )=\lambda f(a)+(1-\lambda )f(b)-f(\lambda a+(1-\lambda )b)-M\frac{%
\lambda (1-\lambda )}{2}(a-b)^{2}.
\end{equation*}%
Since $h(0)=h(1)=0$, 
\begin{equation*}
0=(1-\lambda )h(0)+\lambda h(1)\geq h((1-\lambda )\cdot 0+\lambda \cdot
1)=h(\lambda ),
\end{equation*}%
for all $\lambda \in \lbrack 0,1]$. The assertion is now clear.

\hfill \hbox{\rule{6pt}{6pt}}

For a slight generalization and alternative proof of Theorem \ref{th_2} the
reader is referred to \cite[Theorem 4.2]{FMM2011}.\medskip\ 

\begin{Rem}
By integrating each term of the inequality (\ref{ref_1}) on $\left[ 0,1%
\right] $ with respect to the variable $\lambda $ we recover the inequality (%
\ref{cpn_2}).
\end{Rem}

\begin{Cor}
Preserving the notation of Theorem \ref{th_2}, the following inequalities
hold:%
\begin{eqnarray}
m\frac{(1-2\lambda )^{2}}{8}(a-b)^{2} &\leq &\frac{f(\lambda a+(1-\lambda
)b)+f((1-\lambda )a+\lambda b)}{2}-f\left( \frac{a+b}{2}\right)
\label{ref_3} \\
&\leq &M\frac{(1-2\lambda )^{2}}{8}(a-b)^{2}
\end{eqnarray}%
for all $\lambda \in \lbrack 0,1].$
\end{Cor}

\textit{Proof:} According to Theorem \ref{th_2} for $\lambda =\frac{1}{2}$
we obtain the following result, previously established in \cite{AF2010}: 
\begin{equation}
\frac{m}{8}(b-a)^{2}\leq \frac{f(a)+f(b)}{2}-f\left( \frac{a+b}{2}\right)
\leq \frac{M}{8}(b-a)^{2}.  \label{ref_2}
\end{equation}%
We consider the above inequality (\ref{ref_2}) replacing $a\rightarrow
\lambda a+(1-\lambda )b$ and $b\rightarrow (1-\lambda )a+\lambda b$ (the
hypothesis $m\leq f^{\prime \prime }\leq M$ is still working on the interval
with these endpoints because it is contained by $[a,b]$) and we get the
claimed result.

\hfill \hbox{\rule{6pt}{6pt}}

\begin{Rem}
Notice that by integrating all terms of (\ref{ref_3}) on $\left[ 0,1\right] $
with respect to $\lambda $ we recover now the inequality (\ref{cpn_1}).
\end{Rem}

Next we give some estimates of the Fej\'{e}r inequalities (Proposition \ref%
{prop_fej}):

\begin{The}
\label{fej_2}Let $f:[a,b]\rightarrow \mathbb{R}$ be a twice differentiable
function such that there exist real constants $m$ and $M$ so that $m\leq
f^{\prime \prime }\leq M$. Assume $g:[a,b]\rightarrow \mathbb{R}_{+}$\textit{%
\ }{\ is integrable and symmetric about }$\frac{a+b}{2}.$ Then{\ }the
following inequalities hold:%
\begin{eqnarray}
\frac{m}{2}\int_{a}^{b}(t-a)(b-t)g\left( t\right) \mathrm{d}t &\leq &\frac{%
f\left( a\right) +f\left( b\right) }{2}\int_{a}^{b}g\left( t\right) \mathrm{d%
}t-\int_{a}^{b}f(t)g\left( t\right) \mathrm{d}t\   \label{1} \\
&\leq &\frac{M}{2}\int_{a}^{b}(t-a)(b-t)g\left( t\right) \mathrm{d}t\ \ 
\end{eqnarray}%
and%
\begin{eqnarray}
\frac{m}{8}\int_{a}^{b}(2t-a-b)^{2}g\left( t\right) \mathrm{d}t &\leq
&\int_{a}^{b}f(t)g\left( t\right) \mathrm{d}t-f\left( \frac{a+b}{2}\right) \
\int_{a}^{b}g\left( t\right) \mathrm{d}t  \label{2} \\
&\leq &\frac{M}{8}\int_{a}^{b}(2t-a-b)^{2}g\left( t\right) \mathrm{d}t.\ \ 
\end{eqnarray}
\end{The}

\textit{Proof:} We multiply (\ref{ref_1}) by $g\left( \lambda a+(1-\lambda
)b\right) $ and integrate the result on $\left[ 0,1\right] $ with respect to
the variable $\lambda $. Using the change of the variable $\lambda
a+(1-\lambda )b=t$ we get 
\begin{eqnarray}
&&\frac{m}{2}\int_{a}^{b}(t-a)(b-t)g\left( t\right) \mathrm{d}t  \notag \\
&\leq &f\left( a\right) \int_{a}^{b}\frac{b-t}{b-a}g\left( t\right) \mathrm{d%
}t+f\left( b\right) \int_{a}^{b}\frac{t-a}{b-a}g\left( t\right) \mathrm{d}%
t-\int_{a}^{b}f(t)g\left( t\right) \mathrm{d}t\   \notag \\
&\leq &\frac{M}{2}\int_{a}^{b}(t-a)(b-t)g\left( t\right) dt.
\label{right_fejer_1}
\end{eqnarray}%
On the other hand, due to the symmetry property of $g$, for $t=a+b-x,$ we
also have 
\begin{eqnarray}
&&\frac{m}{2}\int_{a}^{b}(b-x)(x-a)g\left( x\right) dx  \notag \\
&\leq &f\left( a\right) \int_{a}^{b}\frac{x-a}{b-a}g\left( x\right)
dx+f\left( b\right) \int_{a}^{b}\frac{b-x}{b-a}g\left( x\right)
dx-\int_{a}^{b}f(t)g\left( t\right) dt\   \notag \\
&\leq &\frac{M}{2}\int_{a}^{b}(b-x)(x-a)g\left( x\right) dx.
\label{right_fejer_2}
\end{eqnarray}%
Summing (\ref{right_fejer_1}) and (\ref{right_fejer_2}) we find (\ref{1}).

In order to prove the remaining inequalities (\ref{2})\textit{\ }we follow
same steps as above, using (\ref{ref_3}) instead of (\ref{ref_1}). The
computation is straightforward, taking into account the symmetry of $g$
(applied now as $g\left( \lambda a+(1-\lambda )b\right) =g\left( (1-\lambda
)a+\lambda b\right) $). We omit the details.

This completes the proof.

\hfill \hbox{\rule{6pt}{6pt}}

It is remarkable that (\ref{1}) agrees, having an extended form, with \cite[%
pp.53, Exercise 4]{CPN2006}.

\begin{Rem}
If $g:\left[ a,b\right] \rightarrow \left[ 0,1\right] $ then the function $%
h\left( x\right) =1-g\left( x\right) $ satisfies the same symmetry and
positivity conditions and Theorem \ref{fej_2} also applies. That yields the
following estimates of the precision in (\ref{1}): 
\begin{eqnarray}
&&\left( b-a\right) \left( \frac{f\left( a\right) +f\left( b\right) }{2}-%
\frac{1}{b-a}\int_{a}^{b}f(t)\mathrm{d}t\ -m\frac{(b-a)^{2}}{12}\right) 
\notag \\
&\geq &\frac{f\left( a\right) +f\left( b\right) }{2}\int_{a}^{b}g\left(
t\right) \mathrm{d}t-\int_{a}^{b}f(t)g\left( t\right) \mathrm{d}t\ -\ \frac{m%
}{2}\int_{a}^{b}(t-a)(b-t)g\left( t\right) \mathrm{d}t\geq 0
\end{eqnarray}%
and%
\begin{eqnarray}
&&\left( b-a\right) \left( M\frac{(b-a)^{2}}{12}-\frac{f\left( a\right)
+f\left( b\right) }{2}+\frac{1}{b-a}\int_{a}^{b}f(t)\mathrm{d}t\ \right) 
\notag \\
&\geq &\frac{M}{2}\int_{a}^{b}(t-a)(b-t)g\left( t\right) \mathrm{d}t-\frac{%
f\left( a\right) +f\left( b\right) }{2}\int_{a}^{b}g\left( t\right) \mathrm{d%
}t+\int_{a}^{b}f(t)g\left( t\right) \mathrm{d}t\ \geq 0.\ 
\end{eqnarray}%
By a similar technique one can estimate (\ref{2}).
\end{Rem}

\begin{Rem}
For the particular case $g\left( x\right) =1$ if we apply Theorem \ref{fej_2}
on the intervals $\left[ a,\frac{a+b}{2}\right] ,$ $\left[ \frac{a+b}{2},b%
\right] $ we get:%
\begin{equation}
\frac{m\left( b-a\right) ^{2}}{48}\leq \frac{1}{2}\left( \frac{f\left(
a\right) +f\left( b\right) }{2}+f\left( \frac{a+b}{2}\right) \right) -\frac{1%
}{b-a}\int_{a}^{b}f(t)g\left( t\right) \mathrm{d}t\leq \frac{M\left(
b-a\right) ^{2}}{48}\ 
\end{equation}%
and 
\begin{equation}
\frac{m\left( b-a\right) ^{2}}{96}\leq \frac{1}{b-a}\int_{a}^{b}f(t)g\left(
t\right) \mathrm{d}t-\frac{1}{2}\left( f\left( \frac{3a+b}{4}\right) \
+f\left( \frac{a+3b}{4}\right) \right) \leq \frac{m\left( b-a\right) ^{2}}{96%
}.\ \ 
\end{equation}
\end{Rem}

The following theorem gives new Fej\'{e}r-type inequalities.

\begin{The}
\label{th_3}Let $f:[a,b]\rightarrow \mathbb{R}$ be a continuous, convex
function and $g:[a,b]\rightarrow \mathbb{R}_{+}$ be continuous. Then the
following statements hold.

1) If $g$ is monotonically decreasing then%
\begin{equation}
\frac{f\left( a\right) +f\left( b\right) }{2}\int_{a}^{b}g\left( t\right) 
\mathrm{d}t-\int_{a}^{b}f(t)g\left( t\right) \mathrm{d}t\ \geq \frac{f\left(
a\right) +f\left( x\right) }{2}\int_{a}^{x}g\left( t\right) \mathrm{d}%
t-\int_{a}^{x}f(t)g\left( t\right) \mathrm{d}t\ \geq 0;  \label{fg_1}
\end{equation}

2) If $g$ is monotonically increasing then 
\begin{equation}
\int_{a}^{b}f(t)g\left( t\right) \mathrm{d}t\ -f\left( \frac{a+b}{2}\right)
\int_{a}^{b}g\left( t\right) \mathrm{d}t\geq \int_{a}^{x}f(t)g\left(
t\right) \mathrm{d}t\ -f\left( \frac{a+x}{2}\right) \int_{a}^{x}g\left(
t\right) \mathrm{d}t\geq 0  \label{fg_2}
\end{equation}%
for all $x\in \left( a,b\right) .$
\end{The}

\textit{Proof:} 1) We consider the function $h_{1}:[a,b]\rightarrow \mathbb{R%
}$, defined by 
\begin{equation*}
h_{1}(x)=\frac{f\left( a\right) +f\left( x\right) }{2}\int_{a}^{x}g\left(
t\right) \mathrm{d}t-\int_{a}^{x}f(t)g\left( t\right) \mathrm{d}t\ .
\end{equation*}%
Its first derivative is%
\begin{equation*}
h_{1}^{\prime }(x)=\frac{f^{\prime }\left( x\right) }{2}\int_{a}^{x}g\left(
t\right) \mathrm{d}t-\frac{f\left( x\right) -f\left( a\right) }{2}g\left(
x\right) .
\end{equation*}%
Using the mean value theorems there exist $c_{1},c_{2}\in \left[ a,x\right] $
such that%
\begin{equation*}
h_{1}^{\prime }(x)=\left( \frac{f^{\prime }\left( x\right) }{2}g\left(
c_{1}\right) -\frac{f^{\prime }\left( c_{2}\right) }{2}g\left( x\right)
\right) \left( x-a\right) .
\end{equation*}%
Thus, by the convexity of $f$ and to the monotonicity of $g,$ we have $%
f^{\prime }\left( x\right) \geq f^{\prime }\left( c_{2}\right) $ and $%
g\left( c_{1}\right) \geq g\left( x\right) ,$ hence we conclude that $h_{1}$
is increasing on its domain and $h_{1}\left( b\right) \geq h_{1}\left(
x\right) \geq h_{1}\left( a\right) =0.$ Thus we have (\ref{fg_1}), as
asserted.

2) Similarly, we consider the function $h_{2}:[a,b]\rightarrow \mathbb{R}$,
defined by 
\begin{equation*}
h_{2}(x)=\int_{a}^{x}f(t)g\left( t\right) \mathrm{d}t\ -f\left( \frac{a+x}{2}%
\right) \int_{a}^{x}g\left( t\right) \mathrm{d}t
\end{equation*}%
and we compute its first derivative 
\begin{equation*}
h_{2}^{\prime }(x)=\left( f(x)-f\left( \frac{a+x}{2}\right) \right) g\left(
x\right) \ -\frac{1}{2}f^{\prime }\left( \frac{a+x}{2}\right)
\int_{a}^{x}g\left( t\right) \mathrm{d}t.
\end{equation*}%
There exist $k_{1}\in \left[ \frac{a+x}{2},x\right] $ and $k_{2}\in \left[
a,x\right] $ such that%
\begin{equation*}
h_{2}^{\prime }(x)=\left( f^{\prime }\left( k_{1}\right) g\left( x\right) \
-f^{\prime }\left( \frac{a+x}{2}\right) g\left( k_{2}\right) \right) \frac{%
x-a}{2}.
\end{equation*}%
Therefore, due to the monotonicity we have $f^{\prime }\left( k_{1}\right)
\geq f^{\prime }\left( \frac{a+x}{2}\right) $ and $g\left( x\right) \geq
g\left( k_{2}\right) ,$ which leads that $h_{2}$ is increasing on its domain
and $h_{2}\left( b\right) \geq h_{2}\left( x\right) \geq h_{2}\left(
a\right) =0.$

Thus the proof is completed.

\hfill \hbox{\rule{6pt}{6pt}}

The following direct consequence incorporates the classic statement of the
Hermite-Hadamard inequality.

\begin{Cor}
\label{cor}Suppose $f:[a,b]\rightarrow \mathbb{R}$ is continuous and convex.
Then%
\begin{equation}
\frac{f\left( a\right) +f\left( b\right) }{2}-\frac{1}{b-a}\int_{a}^{b}f(t)%
\mathrm{d}t\ \geq \frac{x-a}{b-a}\left( \frac{f\left( a\right) +f\left(
x\right) }{2}-\frac{1}{x-a}\int_{a}^{x}f(t)\mathrm{d}t\ \right) \geq 0
\label{c_1}
\end{equation}%
and 
\begin{equation}
\frac{1}{b-a}\int_{a}^{b}f(t)\mathrm{d}t\ -f\left( \frac{a+b}{2}\right) \geq 
\frac{x-a}{b-a}\left( \frac{1}{x-a}\int_{a}^{x}f(t)\mathrm{d}t\ -f\left( 
\frac{a+x}{2}\right) \right) \geq 0  \label{c_2}
\end{equation}%
for all $x\in \left( a,b\right) .$
\end{Cor}

\textit{Proof:} We apply Theorem \ref{th_3} with $g\left( x\right) =1,$
which satisfies both monotonicity conditions.

\hfill \hbox{\rule{6pt}{6pt}}

\begin{Rem}
Under the same assumptions as in Proposition \ref{nic} when we apply
Corollary \ref{cor} to the convex function $f(x)-m\frac{x^{2}}{2},$ we
recover and improve the inequalities (\ref{cpn_1}) as follows:%
\begin{eqnarray}
&&\frac{1}{b-a}\int_{a}^{b}f(t)\mathrm{d}t-f\left( \frac{a+b}{2}\right) -m%
\frac{(b-a)^{2}}{24}  \notag \\
&\geq &\frac{x-a}{b-a}\left( \frac{1}{x-a}\int_{a}^{x}f(t)\mathrm{d}%
t-f\left( \frac{a+x}{2}\right) -m\frac{(x-a)^{2}}{24}\right) \geq 0
\end{eqnarray}%
and%
\begin{eqnarray}
&&\frac{M}{8}(b-a)^{2}-\left( \frac{1}{b-a}\int_{a}^{b}f(t)\mathrm{d}%
t-f\left( \frac{a+b}{2}\right) \right)  \notag \\
&\geq &\frac{x-a}{b-a}\left[ \frac{M}{8}(x-a)^{2}-\left( \frac{1}{x-a}%
\int_{a}^{x}f(t)\mathrm{d}t-f\left( \frac{a+x}{2}\right) \right) \right]
\geq 0.\ 
\end{eqnarray}%
\textit{\ Similarly if we use the convex function }$M\frac{x^{2}}{2}-f\left(
x\right) $ we get improvements of (\ref{cpn_2}) which at this moment can
easily be written by the interested reader. Obviously same steps could be
followed from Theorem \ref{th_3}, improving that way Theorem \ref{fej_2}.
\end{Rem}

We end this section with the weighted statement of a known result concerning
convex functions.

In the light of Proposition \ref{prop_fej}, the following statement appears
as a trivial generalization of a result due to Vasi\'{c} and Lackovi\'{c} 
\cite{vas76}, and Lupa\c{s} \cite{lup76} (cf. J. E. Pe\v{c}ari\'{c} et al. 
\cite[pp. 143]{pec92}) and we omit its proof.

\begin{Prop}
Let $p$ and $q$ be two positive numbers and $a_{1}\leq a\leq b\leq b_{1}.$
Let $g:[a,b]\rightarrow \mathbb{R}_{+}$\textit{\ }{\ be integrable and
symmetric about }$A=\frac{pa+qb}{p+q}.$ Then the inequalities%
\begin{equation}
f\left( \frac{pa+qb}{p+q}\right) \int_{A-y}^{A+y}g\left( t\right) \mathrm{d}%
t\leq \int_{A-y}^{A+y}\,f(x)\,g\left( x\right) \mathrm{d}x\leq \frac{%
pf\left( a\right) +qf\left( b\right) }{p+q}\int_{A-y}^{A+y}g\left( t\right) 
\mathrm{d}t
\end{equation}%
hold for $y>0$ and all continuous convex functions $f:\left[ a_{1},b_{1}%
\right] \rightarrow \mathbb{R}$ if and only if 
\begin{equation*}
y\leq \frac{b-a}{p+q}\min \left\{ p,q\right\} .
\end{equation*}
\end{Prop}

\section{Application to special means}

From the inequality (\ref{c_1}) applied to the convex function $t^{p}$, with 
$p\in \left( -\infty ,0\right) \cup \left[ 1,\infty \right) \setminus
\left\{ -1\right\} $ we have$\ $%
\begin{equation}
\left( b-a\right) \left\{ \left[ A_{p}\left( a,b\right) \right] ^{p}-\left[
L_{p}\left( a,b\right) \right] ^{p}\right\} \geq \left( x-a\right) \left\{ %
\left[ A_{p}\left( a,x\right) \right] ^{p}-\left[ L_{p}\left( a,x\right) %
\right] ^{p}\right\} ,  \label{AL}
\end{equation}%
where $x\in \left[ a,b\right] .$ Here $A_{p}\left( a,b\right) =\left( \frac{%
a^{p}+b^{p}}{2}\right) ^{1/p}$ is the power mean and $L_{p}\left( a,b\right)
=\left( \frac{b^{p+1}-a^{p+1}}{\left( p+1\right) \left( b-a\right) }\right)
^{1/p}$ is the $p$-logarithmic mean. Also the limit case $p\rightarrow -1$
(or we may equivalently say the case of the convex function $1/t$) gives us%
\begin{equation*}
\left( b-a\right) \left\{ \frac{1}{H\left( a,b\right) }-\frac{1}{L\left(
a,b\right) }\right\} \geq \left( x-a\right) \left\{ \frac{1}{H\left(
a,x\right) }-\frac{1}{L\left( a,x\right) }\right\} ,
\end{equation*}%
where $H\left( a,b\right) =\frac{2ab}{a+b}$ is the harmonic mean and $%
L\left( a,b\right) =\frac{b-a}{\log b-\log a}$ is the logarithmic mean.

It is also useful to consider the inequality (\ref{c_1}) applied for the
convex function $-\log t$, when we get 
\begin{equation*}
\left[ \frac{A\left( a,b\right) }{I\left( a,b\right) }\right] ^{b-a}\geq %
\left[ \frac{A\left( a,x\right) }{I\left( a,x\right) }\right] ^{x-a},
\end{equation*}%
for $a\neq b,$ $x\in \left[ a,b\right] ,$ where $A\left( a,b\right) =\frac{%
a+b}{2}$ is the arithmetic mean and $I\left( a,b\right) =\frac{1}{e}$ $%
\left( \frac{b^{b}}{a^{a}}\right) ^{\frac{1}{b-a}}$ is the identric mean.

In the remainder, we focus on two immediate particular cases of Theorem \ref%
{th_2} that help us to give improvements of the well known
arithmetic-geometric mean inequality (also known as Young's inequality).

\emph{1)} We apply the theorem to the function $f:[a,b]\rightarrow \mathbb{R}%
\ (a>0)$ defined by $f(x)=-\log x$, which leads to 
\begin{equation}
e^{\frac{\lambda (1-\lambda )(a-b)^{2}}{2b^{2}}}\leq \frac{\lambda
a+(1-\lambda )b}{a^{\lambda }b^{1-\lambda }}\leq e^{\frac{\lambda (1-\lambda
)(a-b)^{2}}{2a^{2}}}.  \label{3}
\end{equation}%
Since $e^{\frac{\lambda (1-\lambda )(a-b)^{2}}{2b^{2}}}\geq 1$, we obtain a
refinement of Young's inequality where $\lambda \in \lbrack 0,1].$

We also obtained a reverse inequality for Young's inequality.

\emph{2)} Next, we apply the theorem to the function $f:[\log a,\log
b]\rightarrow \mathbb{R}$, defined by $f(x)=\exp x$ and we arrive at 
\begin{equation}
\frac{\lambda (1-\lambda )a}{2}\log ^{2}\left( \frac{a}{b}\right) \leq
\lambda a+(1-\lambda )b-a^{\lambda }b^{1-\lambda }\leq \frac{\lambda
(1-\lambda )b}{2}\log ^{2}\left( \frac{a}{b}\right) ,  \label{5}
\end{equation}%
where $a,b>0$ and $\lambda \in \lbrack 0,1]$.

The inequality (\ref{5}) gives an improvement of Young's inequality.

\section*{Acknowledgements}

The first author was supported in part by the Romanian Ministry of
Education, Research and Innovation through the PNII Idei project 842/2008.
The second author was supported by CNCSIS Grant $420/2008.$

\noindent

\end{document}